\documentclass[journal]{IEEEtran}
\usepackage{cite}

\ifCLASSINFOpdf
\else
\fi
\usepackage[cmex10]{amsmath}

\usepackage{algorithmic}
\usepackage{algorithm}
\usepackage{url}
\usepackage{epsfig,epsf,psfrag,amssymb,amsfonts,latexsym,graphicx,mathrsfs,float,multirow,tabularx}
\usepackage[dvips]{color}
\usepackage{array}
\usepackage{subfigure}
\newcommand*{\QEDB}{\hfill\ensuremath{\square}}  %

\IEEEoverridecommandlockouts
\begin{document}
%
\title{Coordinated Multi-area Economic Dispatch via Critical Region Projection}
%
%
%

\author{Ye~Guo,~\IEEEmembership{Member,~IEEE,}
        Lang~Tong,~\IEEEmembership{Fellow,~IEEE,}
        Wenchuan~Wu,~\IEEEmembership{Senior~Member,~IEEE,}
        and Boming~Zhang,~\IEEEmembership{Fellow,~IEEE,}
        and~Hongbin~Sun,~\IEEEmembership{Senior~Member,~IEEE}
\thanks{This work was supported in part by the National Science Foundation under Grant 15499. Part of this work was presented at 2016 American Control Conference.}
\thanks{\scriptsize
Ye Guo and Lang Tong are with the School of Electrical and Computer Engineering, Cornell University, Ithaca, NY, USA.(Email: \{yg299,lt35\}@cornell.edu).}
\thanks{\scriptsize Wenchuan Wu, Boming Zhang, and Hongbin Sun are with the Department of Electrical Engineering, Tsinghua University, Beijing, China.}
}

\maketitle

\begin{abstract}
A coordinated economic dispatch method for multi-area power systems is proposed. Choosing boundary phase angles as coupling variables, the proposed method exploits the structure of critical regions in local problems defined by active and inactive constraints. For a fixed boundary state given by the coordinator, local operators compute the coefficients of critical regions containing the boundary state and of the optimal cost functions then communicate them to the coordinator who in turn optimizes the boundary state to minimize the overall cost. By iterating between local operators and the coordinator, the proposed algorithm converges to the global optimal solution in finite steps, and it requires limited information sharing.
\end{abstract}

%
\IEEEpeerreviewmaketitle

\section{Introduction}
\subsection{Motivation}
Large interconnected power systems are often operated by independent system operators (ISOs), each has its own operating area within which internal resources are used economically.  The operating areas are connected physically by tie-lines that allow one area to import from or export to neighboring areas for better utilization of overall system resources.
Existing approaches to tie-line scheduling rely on trades across borders at proxy buses by market participants. The ad hoc uses of proxy buses and the imperfect information used by market participants result in substantial economic loss, estimated at the level of \$784 million annually for the New York and New England customers \cite{White&Pike:11WP}.

Ideally, the optimal utilization of tie-lines is determined by the {\em joint economic dispatch} (JED) that treats interconnected operating areas as one.  Because each operating area is controlled by an ISO, joint optimality  needs to be achieved in a decentralized fashion, possibly involving a coordinator.  Typically, each ISO optimizes its internal dispatch and exchanges intermediate solutions with its neighbors or the coordinator. This process iterates until convergence. One of the major challenges of implementing decentralized (but jointly optimal) economic dispatch is to limit the number of iterations without involving each area discloses its private information \cite{PJMNY_Agreement}.

\subsection{Related Works}
Multi-area economic dispatch (MAED) has been studied extensively, dating  back to \cite{Lin&Viviani:1984TPAS} in the 1980's.   Existing techniques can be classified based on the methodology used in decomposing decision variables.
The primal decomposition methods partition the decision variables of the overall problem into local and coupling variables.  The dual decomposition techniques, on the other hand, solve a relaxed local problem and use the dual variables to coordinate local optimizations.

Among the dual decomposition methods, the most classic kind of approach is based on Lagrangian relaxation
\cite{Kim&Baldick:97TPS,ConejoAguado98TPS,Chen&Thorp&Mount:03HICCS,Binetti&Etal:14TPS,Erseghe:15TPS,Baldick99APPOPF,Wang&Song&Lu:01IEE,Lai&Xie:15TPS}.  These techniques typically require updates on the Lagrange multipliers, which, depending on the parameter setting, often require a large number of iterations and substantial computation and communication costs.

There is also a body of works based on primal decompositions where coupling variables are first fixed in the subproblems and then solved iteratively as part of the master problem \cite{Nogales&Prieto&Conejo:03OR,Bakirtzis&Biskas:03TPS,Min&Abur:06PSCE,Zhao&LitvinovZheng:14TPS,Chatterjee&Baldic:14COR,Li&Wu&Zhang&Wang::15TPS}.
The key step is to define the coupling variables that need to be solved in the master problem.

Among the primal decomposition algorithms, the recent work of Zhao, Litvinov, and Zheng \cite{Zhao&LitvinovZheng:14TPS} has the special property that the algorithm converges in {\em a finite number of steps}, which is especially attractive for the MAED problem. The key idea in \cite{Zhao&LitvinovZheng:14TPS} is the so-called marginal equivalent decomposition (MED) of variables involving  the set of active constraints and ``free variables''  of the local solutions. By communicating these ``marginal variables'', the algorithm implicitly exploits the finiteness of the structure of active constraint set.
\subsection{Summary of contributions}
In this paper, we propose a MAED method referred to as critical region projection (CRP). As a primal decomposition method, CRP defines for each area a sub-problem using the internal generation as decision variables and its boundary phase angles as coupling variables. The proposed approach is based on a key property in multi-parametric quadratic programming: the optimal generation in each area is a piecewise affine function of boundary state, and its associated optimal cost is a piecewise quadratic function of boundary state. This implies that the space of the boundary state can be partitioned into critical regions, within each region the optimal generation and the optimal cost can be characterized succinctly by the affine and quadratic functions.

CRP iterates between the coordinator and regional operators: Given a boundary state, each area solves its sub-problem, derives its optimal cost as a quadratic function of boundary state, and defines the critical region that contains the given boundary state. The coordinator solves the master problem and projects the point of boundary state to a new critical region with strictly lower cost for the next iteration.

CRP shares some of the important features of the MED approach \cite{Zhao&LitvinovZheng:14TPS}, most important being the finite-step convergence.  Our approach does not require any exchange of system information such as shift factors, status of generations, and capacities of internal generators and branches.
Because the number of boundary buses is relatively small, the parameterization proposed in our approach results in the reduced amount of data exchange. CRP does require a coordinator that may complicate practical implementations.

The reminder of this paper is organized as follows. In section II, we present the JED model and decompose it into local sub-problems and the coordinator's problem. The outline of CRP and the solutions to local sub-problems and the master problem are elaborated in section III. In section IV, we establish the optimality and finite-step convergence of CRP and review its computation/communication costs. In section V, CRP is applied in various test systems and its performance is compared with JED and approaches based on Lagrangian relaxation and MED.

\section{Problem Decomposition}
\subsection{Joint Economic Dispatch Model}
\begin{figure}[h]
  \centering
  \includegraphics[width=0.4\textwidth]{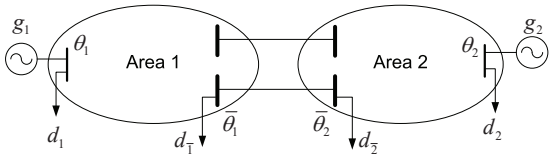}
  \caption{\small An illustration for multi-area power systems}\label{fig:1}
\end{figure}

For simplicity, the MAED model is illustrated via the two area system in Fig.\ref{fig:1}. Similar method can be proposed for systems with more than two areas.

The system state variables are partitioned into four subsets: internal phase angles $\theta_i$ in area $i$ and boundary phase angles $\bar{\theta}_i$ in area $i,(i=1,2)$.

Without loss of generality, we make the following assumptions:

A1) There is no power generation on boundary buses.

A2) Each internal bus is connected with one unit and one load and each boundary bus is connected with one load;

For assumption A1, we can introduce fictitious boundary buses outside the physical ones in case of the presence of boundary generators. With assumption A2, units, loads, and buses have the same indices. Similar approach can be derived if we consider different indices.

The JED is to solve the following centralized optimization:
\begin{align}
&\min \limits_{\{g_{i},\theta_i,\bar{\theta}_i\}}c(g)=\sum_{i=1}^{2}c_{i}(g_{i})=\sum_{i=1}^{2}(g_{i}^T A_{i}g_{i}+b_{i}^T g_{i}),\label{eq:ged_obj}\\
&\textrm{subject}\hspace{0.1cm}\textrm{to}\hspace{0.1cm}H_{i}\theta_{i}+H_{\bar{i}}\bar{\theta}_{i}\leq f_{i},i=1,2,\label{eq:ged_internalline}\\
&\hspace{1.5cm}\bar{H}_{\bar{1}}\bar{\theta}_{1}+\bar{H}_{\bar{2}}\bar{\theta}_{\bar{2}}\leq \bar{f},\label{eq:ged_tieline}\\
&\hspace{1.5cm}\check{g}_{i}\leq g_{i}\leq \hat{g}_{i},i=1,2,\label{eq:ged_ramping}\\
&\left[
\begin{array}{cccc}
Y_{11}&Y_{1\bar{1}}&\mbox{}&\mbox{}\\
Y_{\bar{1}1}&Y_{\bar{1}\bar{1}}&Y_{\bar{1}\bar{2}}&\mbox{}\\
\mbox{}&Y_{\bar{2}\bar{1}}&Y_{\bar{2}\bar{2}}&Y_{\bar{2}2}\\
\mbox{}&\mbox{}&Y_{2\bar{2}}&Y_{22}
\end{array}
\right]\hspace{-0.15cm}\left[
\begin{array}{c}
\theta_{1}\\
\bar{\theta}_{1}\\
\bar{\theta}_{2}\\
\theta_{2}
\end{array}
\right]=\left[
\begin{array}{c}
g_{1}-d_{1}\\
-\bar{d}_{1}\\
-\bar{d}_{2}\\
g_{2}-d_{2}
\end{array}\right],
\label{eq:ged_dclf}
\end{align}
where, as shown in Fig.\ref{fig:1}, the vectors $g_i$ and $d_i$ are internal generations and loads in area $i$ and $\bar{d}_i$ is the vector of boundary load power. The cost functions in (\ref{eq:ged_obj}) are quadratic with coefficients $A_{i}$ and $b_{i}$. The superscript $T$ denotes transpose.

Inequality (\ref{eq:ged_internalline}) represents power flow limits for internal branches of area $i$. Here $H_{i}$ is the branch-bus admittance matrix between internal branches of area $i$ and $\theta_i$, $H_{\bar{i}}$ is the branch-bus admittance matrix between internal branches of area $i$ and $\bar{\theta}_i$, and $f_i$ the power flow limits of internal branches of area $i$. Inequality (\ref{eq:ged_tieline}) describes constraints on boundary power flows, with $\bar{H}_{\bar{i}}$ the branch-bus admittance matrix between tie-lines and boundary state $\bar{\theta}_i$ and $\bar{f}$ the boundary power flow limits. Inequality (\ref{eq:ged_ramping}) restricts the power generations $g_i$ between the lower bound $\check{g}_i$ and upper bound $\hat{g}_i$. Equation (\ref{eq:ged_dclf}) represents the DC load flow equations in which $Y$ is the bus admittance matrix.

In following subsections we decompose the JED model into local optimizations and the coordinator's optimization.
\subsection{Local optimization}
The sub-problem of area $i$ is an economic dispatch (ED) problem with fixed boundary state defined by:
\begin{align}
\min \limits_{\{g_{i},\theta_i\}}&\hspace{0.25cm}c_{i}(g_{i})=g^T_i A_{i}g_{i}+b^T_{i}g_{i},\label{eq:edi_cost}\\
\textrm{subject}\hspace{0.1cm}\textrm{to}&\hspace{0.25cm}H_{i}\theta_{i}+H_{\bar{i}}\bar{\theta}_{i}\leq f_{i},\label{eq:edi_line}\\
\mbox{}&\hspace{0.25cm}\check{g}_i\leq g_{i}\leq \hat{g}_i,\vspace{0.05cm}\label{eq:edi_ramping}\\
\mbox{}&\hspace{-1.0cm}\left[
\begin{array}{cc}
Y_{ii}\hspace{-0.15cm}&\hspace{-0.15cm}Y_{i\bar{i}}\\
Y_{\bar{i}i}\hspace{-0.15cm}&\hspace{-0.15cm}Y_{\bar{i}\bar{i}}
\end{array}\right]\hspace{-0.15cm}\left[\hspace{-0.15cm}
\begin{array}{c}
\theta_{i}\\
\bar{\theta}_{i}
\end{array}\hspace{-0.15cm}\right]\hspace{-0.1cm}=\hspace{-0.1cm}\left[\hspace{-0.1cm}
\begin{array}{c}
g_{i}-d_{i}\\
-\bar{d}_{i}-Y_{\bar{i}\bar{j}}\bar{\theta}_{j}
\end{array}\hspace{-0.1cm}\right].\label{eq:edi_dclf}
\end{align}
By eliminating $\theta_i$ and summarizing all boundary phase angles as $\bar{\theta}=[\bar{\theta}_{i};\bar{\theta}_{j}]$, we write the local sub-problem in area $i$ as
\begin{equation}
\begin{array}{cc}
\min \limits_{g_{i}}&\hspace{0.25cm}c_{i}(g_{i})=g^T_{i}A_{i}g_{i}+b^T_{i}g_{i},\\
\textrm{subject}\hspace{0.1cm}\textrm{to}&\hspace{0.25cm}M_{i}g_{i}+\bar{M}_{i}\bar{\theta}+\tilde{m}_{i}=0,\\
\mbox{}&\hspace{0.25cm}N_{i}g_{i}+\bar{N}_i\bar{\theta}+\tilde{n}_{i}\leq 0,
\end{array}\label{eq:local_MPQP}
\end{equation}
where
\begin{equation}
\hspace{-0.25cm}
\begin{array}{l}
M_{i}=Y_{\bar{i}i}Y_{ii}^{-1},\bar{M}\hspace{-0.1cm}=\hspace{-0.1cm}[Y_{\bar{i}\hspace{0.02cm}\bar{i}}-Y_{\bar{i}i}Y_{ii}^{-1}Y_{i\bar{i}},Y_{\bar{i}\bar{j}}],\\
\tilde{m}_{i}=\bar{d}_{i}-Y_{\bar{i}i}Y_{ii}^{-1}d_{i},
N_{i}=\left[
\begin{array}{c}
H_{i}Y_{ii}^{-1}\\
I\\
-I\end{array}\right],\\
\bar{N}_{i}=\left[\hspace{-0.2cm}
\begin{array}{cc}
-H_{i}Y_{ii}^{-1}Y_{i\bar{i}}\hspace{-0.1cm}+\hspace{-0.1cm}H_{\bar{i}}&\hspace{-0.2cm}0\\
0&\hspace{-0.2cm}0\\
0&\hspace{-0.2cm}0
\end{array}\hspace{-0.2cm}\right],\tilde{n}_{i}\hspace{-0.1cm}=\hspace{-0.1cm}\left[\hspace{-0.2cm}
\begin{array}{c}
-H_{i}Y_{ii}^{-1}d_{i}\hspace{-0.1cm}-\hspace{-0.1cm}f_{i}\\
-\hat{g}_{i}\\
\check{g}_{i}
\end{array}\hspace{-0.2cm}\right].
\end{array}
\end{equation}
Specifically, the equality constraints in (\ref{eq:local_MPQP}) are in the second row of (\ref{eq:edi_dclf}). The inequality constraints in (\ref{eq:local_MPQP}) are arranged in the order of branch power flow limits (\ref{eq:edi_line}) and upper and lower generation limits (\ref{eq:edi_ramping}).

The local sub-problem (\ref{eq:local_MPQP}) has the standard form of multi-parametric quadratic program (MPQP) with boundary phase angles $\bar{\theta}$ as parameters and internal generations $g_i$ as decision variables.

In MPQP, it is of interest to represent the optimal decision variables $g_i^*$ and the value of optimization $c_i(g_i^*)$ as functions of parameters $\bar{\theta}$. Here we give the following theorem that describes the basic properties of the MPQP (\ref{eq:local_MPQP}):

\textit{\textbf{Theorem}} 1  \cite{borrelli2003constrained}: Consider the multi-parametric quadratic programming (\ref{eq:local_MPQP}). Assuming the region $\Theta$ from which the parameters $\bar{\theta}$ take value is convex, then we have the following:

i) The optimal decision variables $g_i^{*}(\bar{\theta})$ is continuous and piecewise affine in $\Theta$;

ii) The value function $J_i^{*}(\bar{\theta})\triangleq c_i(g_i^{*}(\bar{\theta}))$ is continuous, convex, and piecewise quadratic in $\Theta$;

iii) If model (\ref{eq:local_MPQP}) is non-degenerate in $\Theta$, \textit{i.e.}, the rows in matrix $[M_i;\{N_i\}_\mathcal{A}]$ is linearly dependent where $\{N_i\}_\mathcal{A}$ is the sub-matrix of $N_i$ associated with active constraints. Then $J_i^{*}(\bar{\theta})$ is differentiable in $\Theta$.

The key implication of Theorem 1 is that, for the sub-problem of area $i$,  the region $\Theta$ is composed of critical regions. Each critical region corresponds to a particular partition of active and inactive constraints, which is a polyhedron within which $g_i^*(\bar{\theta})$ is an affine function and $J_i^*(\bar{\theta})$ is a quadratic function. Typically, critical regions are half-open-half-closed set. In this paper, to achieve a successive iteration process, we use the closure of critical regions $k$ in the operator $i$'s sub-problem that is denoted as $\Theta_{i,(k)}$. For convenience, we no longer add the word "closure" in the rest of this paper.
\subsection{Coordinator's Optimization}
The main task of the coordinator is to optimize boundary state $\bar{\theta}$ to minimize the overall cost in all areas subjecting to boundary constraints:
\begin{equation}
\begin{array}{ll}
\min \limits_{\bar{\theta}}&J^{*}(\bar{\theta})=\sum \limits_{i=1}^{2}J_{i}^{*}(\bar{\theta}),\\
\textrm{subject}\hspace{0.1cm}\textrm{to}&\bar{H}\bar{\theta}+\tilde{h}\leq 0.
\end{array}\label{eq:11}
\end{equation}
In (\ref{eq:11}) the boundary power flow constraints are written in the same form as local sub-problems in (\ref{eq:local_MPQP}).

The challenge, however, is that the coordinator does not have the exact functional form of $J_i^*$.  Thus (13) cannot be solved directly by the coordinator.  The main idea of CRP, as we describe in the next section,  is to obtain  a partial description of $J_i^*$ from the solution to the local sub-problem $i$ and update boundary state in an iterative fashion.

\section{Proposed Method}
\subsection{Architecture and General Approach}

\begin{figure}[ht]
  \centering
  \includegraphics[width=0.4\textwidth]{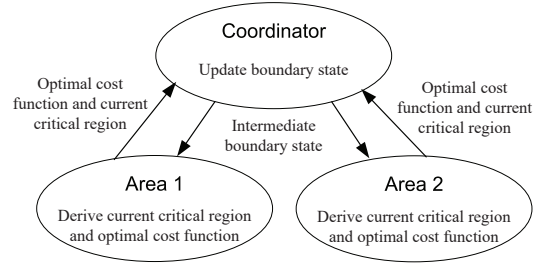}
  \caption{\small The architecture and data flow of CRP}\label{fig:arch}
\end{figure}

We first describe, at a high level, the architecture and the general approach.  As illustrated in Fig.\ref{fig:arch}, the proposed approach involves a coordinator interacting with local area dispatch centers.

Given an intermediate boundary state, each local operator constructs the critical region that contains the boundary state and the parameters of the optimal cost function. Subsequently, the coordinator updates a new boundary state that guarantees a reduced cost for the next iteration.

\begin{figure}[ht]
  \centering
  \includegraphics[width=0.4\textwidth]{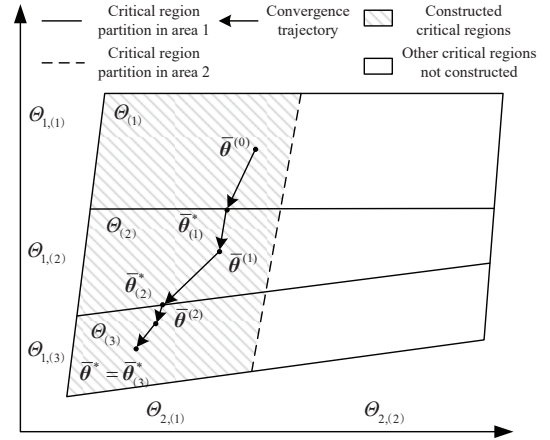}
  \caption{\small Illustration for key steps of CRP}\label{fig:2}
\end{figure}

The detailed constructions of critical regions and projections are described in Sections III.B-C.  Here we illustrate key steps of CRP using a two dimensional example in Fig.\ref{fig:2}.

Initially, the coordinator has the region $\Theta$ from which the boundary state takes value and an initial point $\bar{\theta}^{(0)}\in\Theta$. It communicates $\bar{\theta}^{(0)}$ to areas 1 and 2 who derives the critical regions that contain $\bar{\theta}^{(0)}$, respectively denoted by $\Theta_{1,(1)}$ and $\Theta_{2,(1)}$, and the quadratic optimal cost functions $J_1^*(\bar{\theta})$ and $J_2^*(\bar{\theta})$. The region $\Theta_{(1)}=\Theta_{1,(1)}\cap\Theta_{2,(1)}$ is the critical region of the coordinator's problem in which $J^*(\bar{\theta})$ is quadratic. Hence the coordinator can obtain the optimum point $\bar{\theta}_{(1)}^{*}\in \Theta_{(1)}$ by solving a quadratic programming (QP).

Note that model (\ref{eq:11}) is a convex programming with a unique optimal point. Unless $\bar{\theta}_{(1)}^{*}$ happens to be globally optimal, it resides on the boundary of $\Theta_{(1)}$.

The coordinator then projects the boundary state to a new critical region with strictly lower cost by moving along the anti-gradient direction. See $\bar{\theta}^{(1)}$ in Fig.\ref{fig:2}.  Note that the coordinator does not need the exact form of the new critical region $\Theta_{(2)}$.

In the following iterations, the coordinator sequentially gets $\bar{\theta}^{(2)}$ and $\bar{\theta}^{*}$ along the convergence trajectory shown by the arrows. During the iteration process, we only construct the critical regions through which the convergence trajectory passes, denoted by the shadows. Since there are only finite number of critical regions, the iterative process stops in a finite number of steps.

The following subsections will elaborate the solution to local sub-problems and the method for the coordinator to update the boundary state.

\subsection{Local Sub-problems}\label{sec:CR}
Before elaborating the solution to local sub-problems, we add the following assumptions in CRP:

A3) Given any boundary state that satisfies (\ref{eq:11}), all local sub-problems have feasible solutions;

A4) The JED (\ref{eq:ged_obj})-(\ref{eq:ged_dclf}) has a unique optimal solution;

A5) The local sub-problem (\ref{eq:local_MPQP}) is always non-degenerate.

For assumption A3, the boundary constraints in (\ref{eq:11}) include not only thermal constraints of tie-lines, but also other constraints imposed by system operators (such as limits on maximum export/import power) that guarantee the local sub-problems have feasible solutions. Accordingly, the region $\Theta$ from which the boundary state takes value is defined by
\begin{equation}
\Theta=\{\bar{\theta}|\bar{H}\bar{\theta}+\tilde{h}\leq 0\}. \label{eq:Theta}
\end{equation}

For assumption A5, in case of model (\ref{eq:local_MPQP}) being degenerate, it can be converted to a non-degenerate one by arranging all inequality constraints in a certain sequence, finding as many linearly independent active constraints as possible along the sequence, then setting the other constraints as inactive.

The Lagrangian for the local sub-problem (\ref{eq:local_MPQP}) is
\begin{equation}
\begin{array}{ll}
\!\!\!\!L(g_{i},\lambda_{i},\mu_{i})\!=\!&\!\!\!\!c_{i}(g_i)\!+\!\lambda_{i}^T(M_ig_{i}\!+\!\bar{M}_i\bar{\theta}\!+\tilde{m}_i)\\
\!\!\!\!&\!\!\!+\mu_{i}^{T}(N_ig_i+\bar{N}_i\bar{\theta}+\tilde{n}_i),
\end{array}\label{eq:16}
\end{equation}
where $\lambda_{i}$ and $\mu_{i}$ are the multipliers for the equality and inequality constraints, respectively.
The KKT conditions are
\begin{eqnarray}
\!\!\left[\!\!\!\!
\begin{array}{ccc}
\setlength{\arraycolsep}{0.01pt}
2A_{i}\!\!\!\!&\!\!\!\!M_i^{T}\!\!\!&\!\!\!\{N_i\}_{\mathcal{A}}^{T}\\
M_i\!\!\!\!&\!\!\!\!\mbox{}\!\!\!&\!\!\!\mbox{}\\
\{N_i\}_{\mathcal{A}}\!\!\!\!&\!\!\!\!\mbox{}\!\!\!&\!\!\!\mbox{}
\end{array}
\!\!\!\!\right]\!\!\left[\!\!\!\!
\begin{array}{ccc}
g_{i}\\
\lambda_{i}\\
\{\mu_{i}\}_{\mathcal{A}}
\end{array}\!\!\!\!\right]
\!\!=\!\!\left[\!\!\!
\begin{array}{ccc}
-b_{i}\\
-\bar{M}_i\bar{\theta}-\tilde{m}_i\\
-\{\bar{N}_i\bar{\theta}+\!\tilde{n}_i\}_{\mathcal{A}}
\end{array}
\!\!\!\right]\!\!,
\label{eq:17}\\
\{\mu_{i}\}_{\mathcal{A}}\geq0,\{N_{i}g_{i}+\bar{N}_{i}\bar{\theta}+\tilde{n}_i\}_{\mathcal{A}}=0,\nonumber\\
\{\mu_{i}\}_{\mathcal{I}}=0,\{N_{i}g_{i}+\bar{N}_{i}\bar{\theta}+\tilde{n}_i\}_{\mathcal{I}}\leq0,\nonumber
\end{eqnarray}
where $\{\centerdot\}_{\mathcal{A}}$ and $\{\centerdot\}_{\mathcal{I}}$ denote, respectively, variables associated with active and inactive constraints.

The solution of (\ref{eq:17}) has the form:
\begin{equation}
\left[\!\!
\begin{array}{ccc}
\setlength{\arraycolsep}{0.1pt}
g_{i}\\
\lambda_{i}\\
\{\mu_{i}\}_{\mathcal{A}}
\end{array}\!\!\right]\!\!\!=\!\!\!
\left[\!\!
\begin{array}{ccc}
\setlength{\arraycolsep}{0.1pt}
K_{11}\!\!&\!\!K_{12}\!\!&\!\!K_{13}\\
K_{21}\!\!&\!\!K_{22}\!\!&\!\!K_{23}\\
K_{31}\!\!&\!\!K_{32}\!\!&\!\!K_{33}
\end{array}
\!\!\right]\!\!\!
\left[\!\!\!
\begin{array}{ccc}
-b_{i}\\
-\bar{M}_i\bar{\theta}-\tilde{m}_i\\
-\{\bar{N}_{i}\bar{\theta}+\tilde{n}_i\}_{\mathcal{A}}
\end{array}
\!\!\!\right]\!\!.\label{eq:18}
\end{equation}
\indent For active constraints, their multipliers $\{\mu_{i}\}_{\mathcal{A}}$ are affine functions of $\bar{\theta}$:
\begin{equation}
\begin{array}{ll}
\{\mu_{i}\}_{\mathcal{A}}=&\hspace{-0.25cm}-(K_{32}\bar{M}_i+K_{33}\{\bar{N}_i\}_{\mathcal{A}})\bar{\theta}\\
&\hspace{-0.25cm}-(K_{31}b_i+K_{32}\tilde{m}_i+K_{33}\{\tilde{n}_i\}_{\mathcal{A}})\geq0.
\end{array}
\label{eq:21}
\end{equation}
\indent The optimal generations $g^{*}_{i}$ are also affine functions of $\bar{\theta}$:
\begin{equation}
\begin{array}{l}
g^{*}_{i}=\bar{R}_i\bar{\theta}+\tilde{r}_i,\\
\bar{R}_i=-K_{12}\bar{M}_i-K_{13}\{\bar{N}_{i}\}_{\mathcal{A}},\\
\tilde{r}_i=-K_{11}b_{i}-K_{12}\tilde{m}_i-K_{13}\{\tilde{n}_i\}_{\mathcal{A}}.
\end{array}
\label{eq:22}
\end{equation}
By substituting (\ref{eq:22}) to inactive constraints, we have
\begin{equation}\label{eq:23}
\begin{array}{l}
(\{N_i\}_{\mathcal{I}}\bar{R}_i+\{\bar{N}_i\}_{\mathcal{I}})\bar{\theta}\hspace{-0.05cm}+\hspace{-0.05cm}\{N_i\}_{\mathcal{I}}\tilde{r}_i\hspace{-0.05cm}+\hspace{-0.05cm}\{\tilde{n}_i\}_{\mathcal{I}}\leq 0.
\end{array}
\end{equation}

Given the point of $\bar{\theta}^{(t)}$ and with $g^*_i(\bar{\theta}^{(t)})$, inequality (\ref{eq:21}) defines active constraints via their multipliers, and inequality (\ref{eq:23}) defines inactive constraints via their values.

The intersection of (\ref{eq:21}) and (\ref{eq:23}) defines \emph{current} critical region $k$ that contains $\bar{\theta}^{(t)}$:
\begin{equation}
\begin{array}{l}
\Theta_{i,(k)}=\{\bar{\theta}|\bar{S}_{i,(k)}\bar{\theta}+\tilde{s}_{i,(k)}\leq0\},\\
\bar{S}_{i,(k)}=\left[
\begin{array}{l}
K_{32}\bar{M}_i+K_{33}\{\bar{N}_i\}_{\mathcal{A}}\\
\{N_i\}_{\mathcal{I}}\bar{R}_i+\{\bar{N}_i\}_{\mathcal{I}}
\end{array}
\right],\vspace{0.1cm}\\
\tilde{s}_{i,(k)}=\left[
\begin{array}{l}
K_{31}b_i+K_{32}\tilde{m}_i+K_{33}\{\tilde{n}_i\}_{\mathcal{A}}\\
\{N_i\}_{\mathcal{I}}\tilde{r}_i\hspace{-0.05cm}+\hspace{-0.05cm}\{\tilde{n}_i\}_{\mathcal{I}}
\end{array}
\right].
\end{array}
\label{eq:24}
\end{equation}
The critical region defined by (\ref{eq:24}) is a polyhedron. The redundant inequalities should be removed from (\ref{eq:24}), see \cite{INEQ}.

Within current critical region defined by (\ref{eq:24}), the expression of optimal cost function $J_{i}^{*}(\bar{\theta})$ can be obtained by substituting (\ref{eq:22}) to the cost function (\ref{eq:edi_cost}):
\begin{equation}
J_{i}^{*}(\bar{\theta})=c_i(g_i^{*}(\bar{\theta}))=\bar{\theta}^{T}\bar{A}_{i,(k)}\bar{\theta}+\bar{b}_{i,(k)}^{T}\bar{\theta}+\bar{c}_i,
\label{eq:28}
\end{equation}
where
\begin{equation}
\bar{A}_{i,(k)}=\bar{R}_i^T A_i\bar{R}_i,\bar{b}_{i,(k)}=2\bar{R}_i^{T}A_i\tilde{r}_i+\bar{R}_i^{T}b_i.
\end{equation}
The coordinator knows beforehand that each critical region is a polyhedron and the optimal cost function is quadratic. Therefore, the local system operator only needs to communicate the coefficients $\bar{S}_{i,(k)}$ and $\tilde{s}_{i,(k)}$ in (\ref{eq:24}) and $\bar{A}_{i,(k)}$ and $\bar{b}_{i,(k)}$ in (\ref{eq:28}) to the coordinator.
\subsection{The Coordinator's Problem}
In each iteration, the coordinator searches for the optimal point of $\bar{\theta}$ only within the intersection of current critical regions from local operators:
\begin{align}
\min \limits_{\bar{\theta}}&\hspace{0.2cm}J^{*}(\bar{\theta})=\bar{\theta}^{T}\bar{A}_{\Sigma,(k)}\bar{\theta}+\bar{b}_{\Sigma,(k)}^{T}\bar{\theta}, \label{eq:coobj}\\
\textrm{subject}\hspace{0.1cm}\textrm{to}& \hspace{0.2cm}\bar{S}_{i,(k)}\bar{\theta}+\tilde{s}_{i,(k)}\leq0, \forall i\in 1,2,\label{eq:cocr}\\
\mbox{}& \hspace{0.2cm}\bar{H}\bar{\theta}+\tilde{h}\leq0,\label{eq:cobnd}
\end{align}
where
\begin{equation}
\bar{A}_{\Sigma,(k)}=\sum \limits_{i} \bar{A}_{i,(k)}, \bar{b}_{\Sigma,(k)}=\sum \limits_{i} \bar{b}_{i,(k)}.
\end{equation}
The master problem (\ref{eq:coobj})-(\ref{eq:cobnd}) is a standard QP. CRP converges to the global optimal point $\bar{\theta}^*$ if all constraints associated with critical regions (\ref{eq:cocr}) are inactive. In practise we introduce the stopping tolerance $\epsilon$ on the multipliers $\bar{\mu}$ associated with critical region constraints:
\begin{equation}
\|\bar{\mu}\|_2^2<\epsilon. \label{eq:stop}
\end{equation}

\begin{figure}[htbp]
  \centering
  \includegraphics[width=0.45\textwidth]{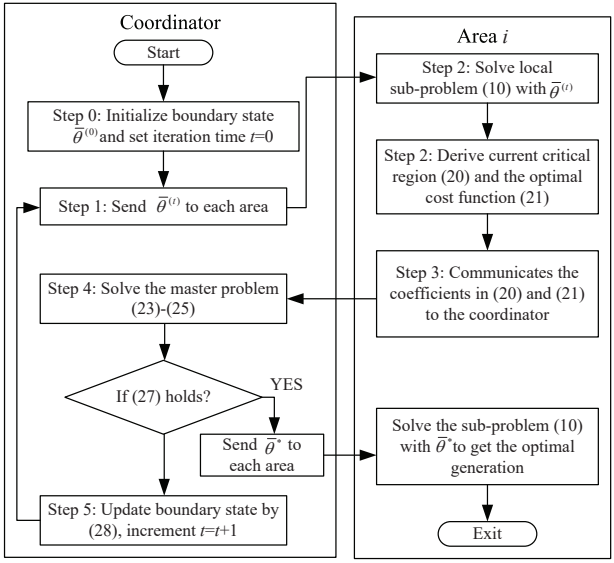}
  \caption{\small The schematic of CRP}\label{fig:scheme}
\end{figure}

If (\ref{eq:stop}) does not hold, then there are active constraints in (\ref{eq:cocr}) and the optimal point in current critical region $k$, denoted by $\bar{\theta}_{(k)}^{*}$, resides on its boundary. According to Theorem 1, the objective function $J^{*}(\bar{\theta})$ is differentiable in $\Theta$. Therefore, the coordinator projects the point of boundary state to a new critical region by moving along the anti-gradient direction:
\begin{equation}
\bar{\theta}^{(t+1)}=\bar{\theta}_{(k)}^{*}-\alpha(P\nabla_{\bar{\theta}}J^{*}),\label{eq:39}
\end{equation}
where $\alpha$ is a small positive constant. The matrix $P$ is the projection matrix that incorporates possible active boundary constraints (\ref{eq:cobnd}), which can be computed by \cite{Haftka_StructOpt}
\begin{equation}
P=I-\{\bar{H}\}_{\mathcal{A}}(\{\bar{H}\}_{\mathcal{A}}^T \{\bar{H}\}_{\mathcal{A}})^{-1} \{\bar{H}\}_{\mathcal{A}}^T.
\end{equation}
The schematic of CRP is given in Fig.\ref{fig:scheme}.
\section{Performance analysis}
We analyze the performance of CRP in this section. In particular, we prove the convergence of CRP and review its
computation/communication costs.
\subsection{Finite-step Convergence and Optimality}\label{sec:opt}
\textit{\textbf{Theorem 2}}: Setting the stopping criterion as (\ref{eq:stop}), we have the following properties on the convergence and optimality of CRP:

i) For any step size $\alpha$ satisfying
\begin{equation}
\alpha < \min \limits_{k}\{\min\{\frac{2}{M_k},l_{k}\}\},\label{eq:upbndalpha}
\end{equation}
where $M_k$ is the maximum eigenvalue of $\bar{A}_{\Sigma,(k)}$ and $l_{k}$ is the distance between $\bar{\theta}_{(k)}^*$ and the boundary of $\Theta$ along the anti-gradient direction at $\bar{\theta}_{(k)}^*$, CRP converges within finite steps, \textit{i.e.}, there exists a constant $K$ such that the iteration of CRP terminates at $t_{\epsilon} <K$;

ii) Assume that the QP solver for the master problem (\ref{eq:coobj})-(\ref{eq:cobnd}) converges to $\epsilon_1$-suboptimality \cite{Boyd_CVXPROG}, \textit{i.e.}, the gap between the objective functions of the primal and dual problems is bounded by

\begin{equation}
J(\bar{\theta}^*_{(k)})-D(\bar{\mu}^*_{(k)},\bar{\nu}^*_{(k)})<\epsilon_1, \label{eq:qpgap}
\end{equation}
where $D$ is the objective function of the dual problem for (\ref{eq:coobj})-(\ref{eq:cobnd}) and $\bar{\nu}$ denotes the multipliers associated with boundary constraints, then the overall cost and generations obtained by CRP converge to the optimal values when $\epsilon$ and $\epsilon_1$ both approach zero, \textit{i.e.},
\begin{equation}
\lim_{\epsilon,\epsilon_1\rightarrow 0} c(g^*(\bar{\theta}^{(t_\epsilon)}))=c(g^*(\bar{\theta}^*)), \label{eq:proof6}
\end{equation}
and
\begin{equation}
\lim_{\epsilon,\epsilon_1\rightarrow 0} g^*(\bar{\theta}^{(t_\epsilon)})=g^*(\bar{\theta}^*). \label{eq:proof15}
\end{equation}

\textit{Proof:}  i) As $J^*$ is convex and piecewise quadratic, consider the entire region of $\Theta$, we have

\begin{equation}
\nabla_{\bar{\theta}}^2 J^* \preceq MI, M=\max \limits_{k} M_k.\label{eq:proof0}
\end{equation}

For (\ref{eq:39}), the values of $J^*(\bar{\theta}_{(k)}^{*})$ and $J^*(\bar{\theta}^{(t+1)})$ yield to
\begin{equation}
\begin{array}{ll}
J^*(\bar{\theta}^{(t+1)})=&\hspace{-0.2cm}J^*(\bar{\theta}_{(k)}^{*})+\nabla_{\bar{\theta}} J^*(\bar{\theta}_{(k)}^{*})^T(-\alpha P\nabla_{\bar{\theta}} J^*(\bar{\theta}_{(k)}^{*}))\\
&\hspace{-0.8cm}+\frac{1}{2}(\alpha P\nabla_{\bar{\theta}} J^*(\bar{\theta}_{(k)}^{*}))^T \nabla_{\bar{\theta}}^2 J^*(z)(\alpha P\nabla_{\bar{\theta}} J^*(\bar{\theta}_{(k)}^{*})).
\end{array}
\label{eq:proof1}
\end{equation}
where $z$ is a point on the line segment between $\bar{\theta}_{(k)}^{*}$ and $\bar{\theta}^{(t+1)}$. To make $J^*(\bar{\theta}^{(t+1)})$ smaller than $J^*(\bar{\theta}_{(k)}^{*})$, the step size $\alpha$ should yield to
\begin{equation}
-\alpha \nabla_{\bar{\theta}} J^*(\bar{\theta}_{(k)}^{*})^T P\nabla_{\bar{\theta}} J^*(\bar{\theta}_{(k)}^{*}) + \alpha^2 \frac{M}{2} \|P\nabla J^*(\bar{\theta}_{(k)}^{*})\|_2^2 < 0. \label{eq:proof11}
\end{equation}
Note that matrix $P$ is idempotent. The solution to (\ref{eq:proof11}) is
\begin{equation}
\alpha < \frac{2}{M}. \label{eq:proof2}
\end{equation}
Furthermore, the point of $J^*(\bar{\theta}^{(t+1)})$ should remain in $\Theta$. Hence the upper bound of step size $\alpha$ is given as (\ref{eq:upbndalpha}). The upper bound in (\ref{eq:upbndalpha}) does not change with iterations.

For any iteration $t$, setting $\alpha$ less than its upper bound, we always have
\begin{equation}
J^*(\bar{\theta}^{(t+1)})<J^*(\bar{\theta}_{(k)}^{*})\leq J^*(\bar{\theta}^{(t)}), \label{eq:proof3}
\end{equation}
which means the objective function strictly decreases by iterations. Furthermore, there are finite number of critical regions and $\bar{\theta}^{(t+1)}$ is in a different critical region from $\bar{\theta}^{(t)}$. Assume that there are $K$ critical regions, then CRP terminates within finite number of iterations $t_{\epsilon}<K$.

ii) The dual problem of the master problem (\ref{eq:coobj})-(\ref{eq:cobnd}) is
\begin{equation}
\begin{array}{l}
\max \limits_{\{\bar{\mu},\bar{\nu}\}} D(\bar{\mu},\bar{\nu})=-\frac{1}{4}([
\bar{H}^T\hspace{0.1cm}\bar{S}_{i,(k)}^T
]\left[
\begin{array}{c}
\bar{\nu}\\
\bar{\mu}
\end{array}\right]+\bar{b}_{\Sigma,(k)})^T\bar{A}_{\Sigma,(k)}^{-1}\\
([
\bar{H}^T\hspace{0.1cm}\bar{S}_{i,(k)}^T
]\left[
\begin{array}{c}
\bar{\nu}\\
\bar{\mu}
\end{array}\right]
+\bar{b}_{\Sigma,(k)})+[
\tilde{h}^T\hspace{0.1cm}\tilde{s}_{i,(k)}^T
]\left[
\begin{array}{c}
\bar{\nu}\\
\bar{\mu}\end{array}\right],\\
\textrm{subject}\hspace{0.1cm}\textrm{to}\hspace{0.1cm}\bar{\nu}\geq0, \bar{\mu}\geq0.
\end{array}\label{eq:dual}
\end{equation}
By substituting (\ref{eq:dual}) to (\ref{eq:qpgap}) we have
\begin{equation}
\begin{array}{l}
J^*(\bar{\theta}^{(t_\epsilon)})-\frac{1}{4}(\bar{\nu}^{(t_\epsilon)})^T\bar{H}\bar{A}_{\Sigma,(k)}^{-1}\bar{H}^T\bar{\nu}^{(t_\epsilon)}\\
-\frac{1}{2}(\bar{\nu}^{(t_\epsilon)})^T\bar{H}\bar{A}_{\Sigma,(k)}^{-1}\bar{S}_{i,(k)}^T\bar{\mu}^{(t_\epsilon)}\\
-\frac{1}{4}(\bar{\mu}^{(t_\epsilon)})^T\bar{S}_{i,(k)}\bar{A}_{\Sigma,(k)}^{-1}\bar{S}_{i,(k)}^T\bar{\mu}^{(t_\epsilon)}\\
-\frac{1}{2}\bar{b}_{\Sigma,(k)}^T\bar{A}_{\Sigma,(k)}^{-1}\bar{H}^T\bar{\nu}^{(t_\epsilon)}
-\frac{1}{2}\bar{b}_{\Sigma,(k)}^T\bar{A}_{\Sigma,(k)}^{-1}\bar{S}_{i,(k)}^T\bar{\mu}^{(t_\epsilon)}\\
-\tilde{h}^T\bar{\nu}^{(t_\epsilon)}-\tilde{s}_{i,(k)}^T\bar{\mu}^{(t_\epsilon)}<\epsilon_1.
\end{array}\label{eq:proof4}
\end{equation}

When CRP terminates at $\bar{\theta}^{(t_\epsilon)}$, by substituting (\ref{eq:stop}) to (\ref{eq:proof4}) and dropping the quadratic term of $\epsilon$, we have

\begin{equation}
\begin{array}{l}
J^*(\bar{\theta}^{(t_\epsilon)}))-\frac{1}{4}(\bar{\nu}^{(t_\epsilon)})^T\bar{H}\bar{A}_{\Sigma,(k)}^{-1}\bar{H}^T\bar{\nu}^{(t_\epsilon)}\\
-\frac{1}{2}\bar{b}_{\Sigma,(k)}^T\bar{A}_{\Sigma,(k)}^{-1}\bar{H}^T\bar{\nu}^{(t_\epsilon)}-\tilde{h}^T\bar{\nu}^{(t_\epsilon)}
<\epsilon_1+\gamma\epsilon.
\end{array}\label{eq:proof5}
\end{equation}
where
\begin{equation}
\begin{array}{ll}
\gamma=&\sup \limits_{\bar{\theta}\in\Theta} (\|\frac{1}{2}(\bar{\nu}^{(t_\epsilon)})^T\bar{H}\bar{A}_{\Sigma,(k)}^{-1}\bar{S}_{i,(k)}^T\\
&+\frac{1}{2}\bar{b}_{\Sigma,(k)}^T\bar{A}_{\Sigma,(k)}^{-1}\bar{S}_{i,(k)}^T+\tilde{s}_{i,(k)}^T\|).
\end{array} \label{eq:gamma}
\end{equation}

Note that inequality (\ref{eq:proof5}) actually bounds the sub-optimality level of the following problem:
\begin{equation}
\begin{array}{ll}
\min \limits_{\bar{\theta}}& J^*(\bar{\theta})=\bar{\theta}^T \bar{A}_{\Sigma,(k)} \bar{\theta}+ \bar{b}_{\Sigma,(k)}^T \bar{\theta}\\
\textrm{subject}\hspace{0.1cm}\textrm{to}& \bar{H}\bar{\theta}+\tilde{h}\leq0.
\end{array}\label{eq:proof9}
\end{equation}
Model (\ref{eq:proof9}) minimizes the overall cost in $\Theta$ by assuming the quadratic function in critical region $k$ holds in the entire region $\Theta$. Let $J'$ be the optimal value for (\ref{eq:proof9}), then from (\ref{eq:proof5}) we have
\begin{equation}
J^*(\bar{\theta}^{(t_\epsilon)})-J' < \epsilon_1+\gamma\epsilon. \label{eq:proof10}
\end{equation}
According to the convexity of $J^*(\bar{\theta})$, there is $J'\leq J^*(\bar{\theta}^*)$. Hence the difference between $J^*(\bar{\theta}^{(t_\epsilon)}))$ and $J^*(\bar{\theta}^*)$ is bounded by
\begin{equation}
J^*(\bar{\theta}^{(t_\epsilon)})-J^*(\bar{\theta}^*) < \epsilon_1+\gamma\epsilon. \label{eq:proof8}
\end{equation}
When $\epsilon$ and $\epsilon_1$ both approach to zero, the limit of the right hand side in (\ref{eq:proof8}) equals to zero. Therefore we have
\begin{equation}
\lim_{\epsilon,\epsilon_1\rightarrow 0} [J^*(\bar{\theta}^{(t_\epsilon)})-J^*(\bar{\theta}^*)]=0. \label{eq:proof13}
\end{equation}
According to the definition of $J^*$, (\ref{eq:proof6}) can be proved. Consequently, (\ref{eq:proof15}) also holds due to the convexity of $c(g)$.
\QEDB

Theorem 2 theoretically proves the convergence and optimality of CRP. In practise, however, we choose $\alpha$ as a small constant according to our experience. We do not really calculate the upper bound in (\ref{eq:upbndalpha}) or the constant $\gamma$ in (\ref{eq:gamma}).

\subsection{Computation/Communication Costs}
The computation cost of CRP mainly includes the following two parts:

i) \emph{Local sub-problem solution and critical region determination in each area}. The local sub-problems have standard forms of QP. The definitions of current critical regions can also be naturally obtained via (\ref{eq:21})-(\ref{eq:24});

ii) \emph{The solution to the master problem (\ref{eq:coobj})-(\ref{eq:cobnd}) at the coordinator}. The master problem also has the standard form of QP. Because the dimension of the QP is the size of the boundary state vector, the computation cost of this step is expected to be small.

On communication cost, as shown in Fig.\ref{fig:arch}, the data exchange in CRP includes the following two parts:

i) \emph{Communications from local areas to the coordinator}.  Each area communicates the coefficients $\bar{S}_{i,(k)}$ and $\tilde{s}_{i,(k)}$ in (\ref{eq:24}) and $\bar{A}_{i,(k)}$ and $\bar{b}_{i,(k)}$ in (\ref{eq:28}) to the coordinator. The number of columns of $\bar{S}_{i,(k)}$ is small, but the numbers of rows of $\bar{S}_{i,(k)}$ and $\tilde{s}_{i,(k)}$ may be large. According to our experience, however, a large portion of the inequalities in (\ref{eq:24}) are redundant and can be eliminated. The sizes of $\bar{A}_{i,(k)}$ and $\bar{b}_{i,(k)}$ equal to the number of boundary buses. In particular, these coefficients do not include any specific information of physical systems.

ii) \emph{Communications from coordinator to local areas.}  The coordinator sends the newest boundary state $\bar{\theta}$ to corresponding areas. This step only involves vector communication.

Furthermore, the finite-step convergence of CRP also guarantees its computation and communication efficiencies.

\section{Numerical Tests}
\subsection{2-area 6-bus system test}

CRP was tested on various test beds and compared with the following three approaches:

i) Direct solution to the JED (\ref{eq:ged_obj})-(\ref{eq:ged_dclf});

ii) The Lagrangian relaxation method (LR) \cite{ConejoAguado98TPS}, the multipliers associated with boundary constraints were
 initialized as zero and the artificial parameters were tuned to achieve relatively fast convergence;

iii) The marginal equivalence decomposition based method (MED) \cite{Zhao&LitvinovZheng:14TPS}, the binding constraints set
 were initialized as void. The quadratic cost functions were approximated by piecewise linear functions with 20 equal size blocks.

In all tests, the initial boundary phase angles of CRP were set as zero. The values for $\epsilon$ and $\epsilon_1$ were set as
$10^{-6}$ and the step size $\alpha$ was set as $10^{-4}$.

\begin{figure}[H]
  \centering
  \includegraphics[width=0.4\textwidth]{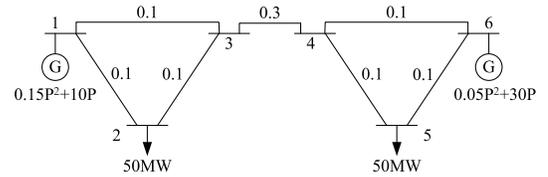}
  \caption{\small Configuration and parameters of 6-bus system}
  \label{fig:6bussys}
\end{figure}

We first compared these four methods on a simple 6-bus system whose configuration, branch reactance, and cost functions were
given in Fig.\ref{fig:6bussys}. The overall costs, iteration times, and computation and communication costs of the four approaches were
compared in TABLE \ref{table:6bussysresults}.
\vspace{-0.7cm}
\begin{table}[H]
\centering
\caption{\small Performances comparison for 6-bus system test}
\begin{tabular}{ccccc}
\hline
Method&Iteration&Overall costs&CPU time&Float data\\
&times&(\$/hr)&costs (ms)&exchanged\\
\hline
JED&-&2375.00&84.28&-\\
LR&12&2376.10&340.92&48\\
MED&2&2375.00&149.27&80*\\
CRP&1&2375.00&113.38&38\\
\hline
\end{tabular}\\
\vspace{0.05cm}
*Shift matrices were not counted, same for other tests
\label{table:6bussysresults}
\end{table}
\vspace{-0.2cm}
LR converged in 12 iterations, its cost was a little higher than that of JED due to the convergence tolerance and its CPU time
 cost was about four times of that of JED. MED needed two iterations to converge to the optimal block in its piecewise linear
 cost functions; its results were optimal in this test and its computation time cost was much less than LR, while its
 communication cost was higher.

On the other hand, as no constraint was considered in this test, there was only one critical region that covered the entire
 boundary state space. Accordingly, CRP achieved the optimal solution within only one iteration. It also had satisfactory
 computation and communication efficiencies.
\vspace{-0.2cm}
\subsection{2-area 44-bus system test}
\begin{figure}[h]
  \centering
  \includegraphics[width=0.41\textwidth]{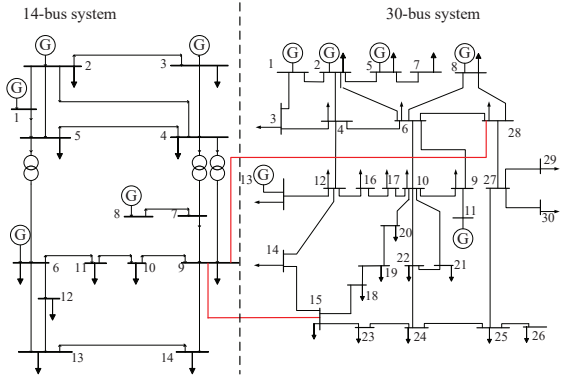}
  \caption{\small Configuration of 14- and 30-bus system}
  \label{fig:sys}
\end{figure}

Similar test was performed on a two area system composed by the IEEE 14- (area 1) and 30-bus (area 2) systems. Two tie-lines were
added between the two areas, the first connected bus 9 in area 1 and bus 15 in area 2 with reactance 0.15p.u., the second
connected bus 9 in area 1 and bus 28 in area 2 with reactance 0.25p.u.. The configuration of the test system was illustrated
in Fig.\ref{fig:sys}:
There were three boundary buses, setting bus 9 in area 1 as phase angle reference, then the space of boundary state had the
dimension of two. The boundary constraints (\ref{eq:ged_tieline}) were
\begin{equation}
\begin{array}{ll}
-50MW\leq P_{9-15}, P_{9-28}\leq 80MW,&\\
-80MW\leq P_{9-15}+P_{9-28}\leq 80MW.&
\end{array}
\end{equation}

Note that IEEE 14- and 30-bus systems are primarily independent and their cost coefficients are very different. Hence
two different scenarios were designed in this test:

i) The cost coefficients in IEEE 30-bus system increased to ten times of their default values.

ii) Default cost coefficients were used.

For both scenarios, the performances of the four approaches were compared in TABLE \ref{table:IEEEresults2}.

\begin{figure}[h]
  \centering
  \includegraphics[width=0.4\textwidth]{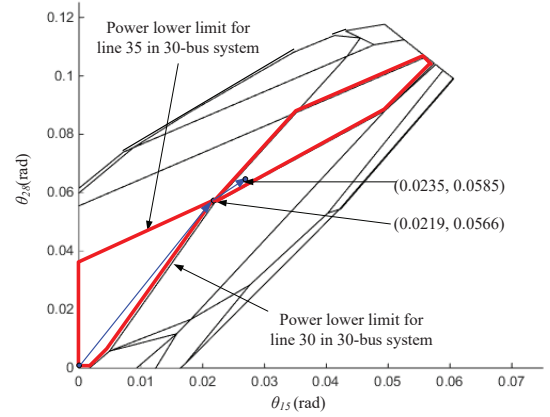}
  \caption{\small The convergence trajectory of CRP in scenario 1}
  \label{fig:add}
\end{figure}

In the first scenario, the prices in the two areas were comparable. Accordingly, the optimum point of boundary state resided
inside $\Theta$ with zero gradient. The CRP method needed two iterations to converge, with one projection of
critical regions. The critical region partition for the boundary state space at the coordinator and the convergence
trajectory were plotted in Fig.\ref{fig:add}. For comparison, LR needed 127 iterations to converge with prohibitive
computation and communication costs. The MED approach converged in three iterations, its results in this test were
sub-optimal due to the piecewise linearization to cost functions. Its CPU time cost was about three times
of that of JED and its communication cost was lower than LR. CRP was the only one out of the three distributed approaches
that achieved the same results with JED, it also needed the least number of iterations and computation/communication costs.
\begin{table}[h]
\centering
\caption{\small Performances comparison for the IEEE system test}
\vspace{-0.3cm}
\vspace{0.3cm}
\begin{tabular}{ccccc}
\hline
Method&Iteration&Overall costs&CPU time&Float data\\
&times&(\$/hr)&costs (ms)&exchanged\\
\hline
&&Scenario 1&\\
\hline
JED&-&14597.54&124.34&-\\
LR&127&14598.11&8933.6&1016\\
MED&3&14599.73&399.43&876\\
CRP&2&14597.54&177.63&188\\
\hline
&&Scenario 2&\\
\hline
JED&-&6095.31&142.74&-\\
LR&270&6095.88&12033.5&2160\\
MED&Infeasible&-&-&-\\
CRP&2&6095.31&183.12&188\\
\hline
\end{tabular}\label{table:IEEEresults2}
\end{table}

\begin{figure}[ht]
  \centering
  \includegraphics[width=0.4\textwidth]{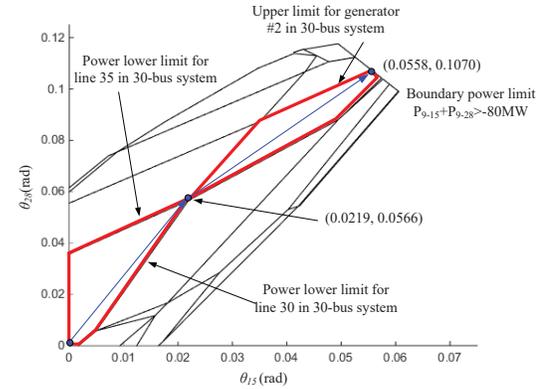}
  \caption{\small The convergence trajectory in scenario 2}
  \label{fig:3a}
\end{figure}

In the second scenario, the prices in area 2 was much lower than those in area 1 and the optimal point of boundary state
resided on the boundary of $\Theta$. The critical region partition and the convergence
trajectory of CRP were given in Fig.\ref{fig:3a}. CRP method needed two iterations to
obtain the same results as JED with reasonable computation and communication costs. For comparison, 
LR needed more iteration times than the first scenario. In MED,
the sub-problem of area 2 became infeasible during its iteration process.
\subsection{3-area 448-bus system test}
\begin{figure}[ht]
  \centering
  \includegraphics[width=0.4\textwidth]{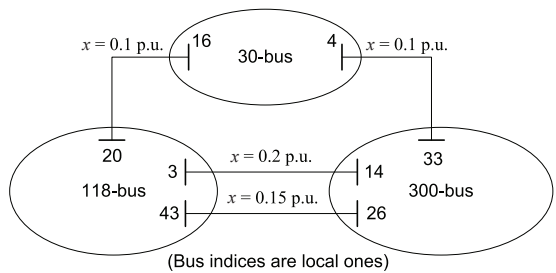}
  \caption{\small Configuration of the 3-area system}
  \label{fig:3area}
\end{figure}

The four MAED approaches were also compared on a 3-area system composed by IEEE 30-bus, 118-bus, and 300-bus systems. Their
interconnections were illustrated in Fig.\ref{fig:3area}. The power limits for all tie-lines were set as 40MW. The performances
of the four approaches were compared in TABLE \ref{table:3areasysresults}.

\begin{table}[h]
\centering
\caption{\small Performances comparison for 3-area 448-bus system test}
\begin{tabular}{ccccc}
\hline
Method&Iteration&Overall costs&CPU time&Float data\\
&times&(\$/hr)&costs (ms)&exchanged\\
\hline
JED&-&$8.31\times10^5$&918.87&-\\
LR&Did not converge&-&-&-\\
MED&5&$8.40\times10^5$&9158.90&3630\\
CRP&5&$8.31\times10^5$&5185.98&1618\\
\hline
\end{tabular}\label{table:3areasysresults}
\end{table}

LR did not converge in this scenario. Both MED and CRP needed five iterations to converge. The overall cost of MED was a
little higher than JED due to the linearization. While CRP got the same cost with JED and needed less computation and
communication costs than MED.

In particular, in the first iteration of CRP, we compared the number of rows in matrices $\bar{S}_{i,(k)}$ before
and after the removal of redundant inequalities as TABLE \ref{table:rcremoval}:

\begin{table}[h]
\centering
\caption{\small The number of inequalities describing current critical regions before and after the redundant removal}
\begin{tabular}{ccc}
\hline
Area&Before&After\\
\hline
30-bus&98&11\\
118-bus&486&19\\
300-bus&966&17\\
\hline
\end{tabular}\label{table:rcremoval}
\end{table}

From TABLE \ref{table:rcremoval} we found that most constraints were redundant. Although CRP might require substantial
communication cost in the worst case, it had satisfactory communication efficiencies in all our simulations. Intuitively,
this is because a low dimensional (the number of boundary buses) polyhedron usually has limited number of edges
(the number of non-redundant constraints).

\subsection{Discussions}

Among the benchmark techniques compared, both CRP and MED require minimum iterations among local operators.
This is a very important feature as the size of local optimization is quite large and the cost of
optimization is substantial. In this respect, the LR technique is at a disadvantage.

Both LR and CRP require minimal information exchange per-iteration. This is also very important in practice.
The MED technique, however, requires local operators to share system parameters and configurations. CRP, on the other hand,
exchange only intermediate boundary state, critical regions and optimal cost functions, which tend to be in low dimensions 
and do not contain any information of internal parts of subareas.

The computation cost of LR per iteration is quite low (although more iterations are needed).
MED and CRP have comparable computation cost, with MED requiring to solve local problems with larger scales and CRP requiring 
computation to obtain critical regions and optimal cost functions.

In summary, experience from our numerical experiments suggested that CRP is competitive in its overall performance 
in accuracy and cost.

\section{Conclusion}
\indent A coordinated multi-area economic dispatch method based on critical region projection is proposed in this paper. With a given boundary state, each area solves its local dispatch problem, determines its current critical region, and derives its optimal cost function. The coordinator minimizes the overall cost within current critical region and then project the boundary state to a new critical region with a reduced cost. The iterative process between local sub-problems and the coordinator will converge to the global optimum solution within finite number of iterations.
{
\bibliographystyle{IEEEtran}
\bibliography{ref}
}

%

%
%
%




\end{document}